\documentclass[12pt]{amsart}
\newtheorem{thm}{Theorem}[section]

\newtheorem{lem}[thm]{Lemma}
\newtheorem{prop}[thm]{Proposition}
\theoremstyle{definition}
\newtheorem{defn}[thm]{Definition}
\theoremstyle{remark}
\newtheorem{rem}[thm]{Remark}
\numberwithin{equation}{section}

\newcommand{\Real}{\mathbb R}
\newcommand{\Sphere}{\mathbb S}

\begin{document}
\centerline{FRANCESCA ANTOCI}
\bigskip
\bigskip
\bigskip
\centerline{{\bf On the spectrum of the Laplace-Beltrami operator  } }
\centerline{{\bf for $p$-forms on asymptotically hyperbolic manifolds}}
\bigskip\bigskip
\date{}
\noindent{\small SUMMARY.-Under suitable conditions on the asymptotic decay of the
metric, we compute the essential spectrum of the Laplace-Beltrami operator
acting on $p$-forms on asymptotically hyperbolic manifolds.
}
\par \bigskip \bigskip \bigskip
\centerline{\small{\bf Sullo spettro dell'operatore di Laplace-Beltrami}}
\centerline{\small{\bf per le $p$-forme su variet\'a asintoticamente iperboliche}}
\par \bigskip
\noindent {\small RIASSUNTO.-
Sotto opportune ipotesi sull'andamento asintotico della metrica, si calcola lo spettro essenziale
dell'operatore di Laplace-Beltrami per le $p$-forme su variet\'a asintoticamente iperboliche.}

\section{Introduction}
The analysis of the spectrum of the Laplace-Beltrami operator on complete
noncompact Riemannian manifolds in its relationships with the geometric
properties of the manifold has been investigated by many
authors. In the case of a general Riemannian manifold the problem turns
out to be very difficult, because of the lack of powerful analytic tools
such as the Fourier transform. Hence the attention has mainly focused on
particular classes of Riemannian manifolds, in which these difficulties
can be bypassed thanks to the presence of symmetries or to the imposition
of a ``controlled" asymptotic behaviour of the Riemannian metric.\par
This is the case for manifolds endowed with rotationally symmetric
Riemannian metrics, where a decomposition technique introduced by Dodziuk
in \cite{Dodziuk} and then employed by Eichhorn
(\cite{Eichhorn}) and Donnelly (\cite{Donnelly}) considerably simplifies
the problem. By this technique, Dodziuk obtained in
\cite{Dodziuk} results on the
existence and multiplicity of $L^2$ harmonic forms for a Riemannian metric
which can be expressed, in geodesic coordinates, as
\begin{equation}\label{metric00}dt^2+ g(t) d\theta^2, \end{equation}
where $g(t)$ is a positive function and $d\theta^2$ is the standard metric
on the sphere $\Sphere^{N-1}$. These techniques were then employed by
 Eichhorn in \cite{Eichhorn} for his results on the discreteness of
the spectrum of the Laplace-Beltrami operator for Riemannian metrics of
type (\ref{metric00}), and by
Donnelly in \cite{Donnelly} in his computation of the spectrum of the
Laplace-Beltrami operator on the hyperbolic space ${\mathbb H}^n$.\par
A completely different approach to this kind of problems can be
found in \cite{Mazzeo}, \cite{Melrose-Mazzeo} and \cite{Mazzeo2}, where
the essential
spectrum is determined on conformally compact Riemannian manifolds
through the sophisticated machinery of the pseudodifferential calculus
on manifolds developed by Melrose (see \cite{Melrose} and the
references therein).\par
In the present paper we consider a noncompact Riemannian
$N$-di\-mensional manifold endowed
with a Riemannian metric of type
\begin{equation}\label{metric01} ds^2=f(t)dt^2 + g(t) d\theta^2,
\end{equation}
where $t\in [0,+\infty)$, $d\theta^2$ is the standard metric on
$\Sphere^{N-1}$, $f(t)>0$ and $g(t)>0$. We suppose that $ds^2$ is
asymptotically hyperbolic, that is $f(t)\rightarrow 1$ and $g(t)
\rightarrow \sinh^2 t$ as $t\rightarrow +\infty$. As for the
behaviour at $t=0$, we suppose that $f(t)=1$ and $g(t)=t^2$ in a
neighbourhood of $0$. Via decomposition and perturbation
techniques, we compute the essential spectrum of the
Laplace-Beltrami operator on $p$-forms, under suitable hypothesis
on the rate of convergence of the metric (\ref{metric01}) to the
hyperbolic metric $$dt^2+\sinh^2t\,d\theta^2.$$ The main result is
the following (Theorem \ref{main}). Let us define $$\tilde f(t):=
f(t)-1,  $$ $$\tilde g(t):=g(t)-\sinh^2 t; $$ if for $t>>0$
\begin{equation}\label{asymestg0}|\tilde g(t)|\leq \frac{C}{t},\quad
|\frac{\partial
\tilde g}{\partial t}|\leq \frac{C}{t},\quad |\frac{\partial^2\tilde
g}{\partial t^2}|\leq \frac{C}{t},\end{equation}
\begin{equation}\label{asymestf0}|\tilde f(t)|\leq
\frac{C}{t}, \quad |\frac{\partial \tilde f}{\partial t}|\leq
\frac{C}{t},\quad|\frac{\partial^2\tilde f}{\partial t^2}|\leq
\frac{C}{t},\end{equation} then the essential spectrum of the
Laplace-Beltrami operator is the interval
$$\left[\min\left\{\left(\frac{N-2p-1}{2}\right)^2,\left(
\frac{N-2p+1}{2}\right)^2\right\},
+\infty\right)$$ if $N\not= 2p$, whilst for $N=2p$ it is equal to
$$\left\{0\right\}\cup \left[\frac{1}{4},+\infty\right).$$
\par We have chosen to express the perturbed metric in form
(\ref{metric01}) in order to have two independent terms in the
perturbation which can be controlled in a clear way.
\par The assumptions (\ref{asymestf0}), (\ref{asymestg0}), though rather general,
can be probably weakened. It would be interesting to get to a more
precise knowledge of the spectrum of $\Delta_M$, in particular as
concerns the absolutely continuous spectrum; however, this seems
difficult, because of the lack of a completely developed Fourier
theory for $p$-forms on the hyperbolic space ${\mathbb H}^N$,
which would permit to understand whether a perturbation of the
Laplace-Beltrami operator is trace-class or not.
\par The paper is organized as follows. In section 2, we construct
an explicit model of asymptotically hyperbolic manifold, endowing
the interior of the unit ball $B^N$ in $\Real^N$ with a Riemannian
metric of type (\ref{metric01}), where $t={\rm
settanh}(\|\overline x\|)$. Moreover, we introduce notations and
some preliminaries which will be useful in the subsequent
sections. In section 3 we prove a generalization of the result by
Dodziuk in \cite{Dodziuk} to the case of a metric of type
(\ref{metric01}); slightly modifying Dodziuk's proof we give
necessary and sufficient conditions for the existence of $L^2$
harmonic $p$-forms on $M$, and we determine their multiplicity. We
then apply the result to the present situation, proving that for
an asymptotically hyperbolic Riemannian manifold $0\in
\sigma_p(\Delta_M)$ if and only if $p=\frac{N}{2}$. Moreover, we
show that in this case $0$ belongs also to the essential spectrum
since it is an eigenvalue of infinite multiplicity. In section 4,
we first introduce an orthogonal decomposition of $L^2_p(M)$
analogous to those employed by Eichhorn and by Donnelly (see
\cite{Eichhorn} and \cite{Donnelly}). The decomposition is
obtained in two steps; first, thanks to the Hodge decomposition on
$\Sphere^{N-1}$, we write any $p$-form $\omega$ as $$\omega=
\omega_{1\delta}\oplus \omega_{2d}\wedge dt\oplus
(\omega_{1d}\oplus \omega_{2\delta}\wedge dt) ,$$ where
$\omega_{1\delta}$ (resp. $\omega_{1d}$) is a coclosed (resp.
closed) $p$-form on $\Sphere^{N-1}$  para\-metrized by $t$, and
$\omega_{2\delta}$ (resp. $\omega_{2d}$) is a coclosed (resp.
closed) $(p-1)$-form on $\Sphere^{N-1}$ parametrized by $t$. The
decomposition is orthogonal in $L^2$ and $\Delta_M$ splits
accordingly as $$\Delta_M=\Delta_{M1}\oplus \Delta_{M2}\oplus
\Delta_{M3}. $$ This allows to reduce ourselves to the study of
the spectral properties of $\Delta_{Mi}$, $i=1,2,3$.\par The
second step consists in decomposing $\omega_{1\delta}$ (resp.
$\omega_{2d}$, $\omega_{2\delta}$) according to an orthonormal
basis of coclosed $p$-eigenforms (resp. closed $(p-1)$-eigenforms,
coclosed $(p-1)$-eigenforms) of $\Delta_{\Sphere^{N-1}}$. In this
way, up to a unitary equivalence, the spectral analysis of
$\Delta_{Mi}$, $i=1,2,3$, can be reduced to the investigation of
the spectra of a countable number of Sturm-Liouville operators
$D_{i\lambda}$ on the half line , parametrized by the eigenvalues
$\lambda$ of $\Delta_{\Sphere^{N-1}}$.\par
In \cite{Eichhorn}
J. Eichhorn proved that for a complete Riemannian metric over a
noncompact manifold the essential spectrum of $\Delta_M$ coincides
with the essential spectrum of the Friedrichs extension
$\Delta_M^F$ of the restriction of $\Delta_M$ to any exterior
domain in $M$. This allows to consider the Sturm-Liouville
operators $D_{i\lambda}$ on $[c,+\infty)$, for $c>0$, and to
overcome the difficulties due to the presence of singular
potentials at $t=0$.\par
In section 5,  under the assumptions (\ref{asymestg0}),
(\ref{asymestf0}), we compute the essential spectrum of
$\Delta_M$. First, through classical perturbation theory, we
compute the spectrum of $D^F_{1\lambda}$ for every
$\lambda$, and we show that
$$\left[\left(\frac{N-2p-1}{2}\right)^2,+\infty\right)\subseteq
\sigma_{\rm ess}(\Delta_{M1}).$$ Then we show that $\sigma_{\rm
ess}(\Delta_{M1})$ is exactly the interval
$\left[\left(\frac{N-2p-1}{2}\right)^2,+\infty\right)$. By duality, we
find that
$\sigma_{\rm ess}(\Delta_{M2})= \left[ \left(\frac{N-2p+1}{2}\right)^2,
+\infty
\right)$. As for the essential spectrum of $\Delta_{M3}$, first we compute
the essential spectrum of $D^F_{3\lambda}$ for every $\lambda$, proving
that $$ \left[\min \left\{\left(\frac{N-2p-1}{2}\right)^2,
\left(\frac{N-2p+1}{2}\right)^2
\right\},+\infty\right) \subseteq \sigma_{\rm
ess}(\Delta_{M3}).$$ Finally we show that any positive number $\mu$ such
that
$$ \mu < \min
\left\{\left(\frac{N-2p-1}{2}\right)^2,\left(\frac{N-2p+1}{2}\right)^2
\right\},$$
can not belong to the essential spectrum of $\Delta_{M3}$. Hence,
$$ \sigma_{\rm ess}(\Delta_M)\setminus \left\{0\right\}=
\left[\min\left\{\left(\frac{N-2p-1}{2}\right)^2,
\left(\frac{N-2p+1}{2}\right)^2\right\},+\infty\right).$$
Then, recalling the results of Section 3, we can fully determine the
essential spectrum of $\Delta_M$.

\section{Preliminary facts}
For $N\geq 2$, let $\overline{B^N}$ denote the closed unit ball
$$\overline{B^N}=\left\{ \bar x= (x_1,...,x_N) \in \Real^N\,|\,
x_1^2+...+x_N^2\leq 1\right\}, $$ and let $\Sphere^{N-1}$ denote
the sphere $$\Sphere^{N-1}= \left\{(x_1,...,x_N)\in \Real^N\,|\,
x_1^2+...+x_N^2=1 \right\} ,$$ endowed with a coordinate system
$(U_i, \Theta_i)$, $i=2,...,k+1$, $\Theta_i: U_i \rightarrow
\Real^{N-1}$.\par Let us consider the interior of
$\overline{B^N}$, $$B^N=\left\{(x_1,...,x_N)\in \Real^N\,|\,
x_1^2+...+x_N^2<1 \right\}, $$ with the coordinate system
$(V_i,\Phi_i)$, for $i=1,...,k+1$, defined in the following way:
in a neighbourhood of $0$, for some $\delta >0$, $$ V_1= \left\{
(x_1,...,x_N)\in \Real^N\,|\, x_1^2 +...+x_N^2 < \delta\right\}$$
and $$ \Phi_1(x_1,...,x_N)=(x_1,...,x_N),$$ whilst for $i>1$,
$\bar x \not=0$, $$V_i=\left\{ \bar x \in \Real^N \,|\, \frac{\bar
x}{\|\bar x\|}\in U_i\right\},$$ $$ \Phi_i: V_i \longrightarrow
(0,+\infty) \times \Theta_i(U_i),$$ $$\Phi_i(x_1,...,x_N)=
\left(2\, {\rm settanh} (\|\bar x\|), \Theta_i\left(\frac{\bar
x}{\|\bar x\|}\right)\right)=:(t,\theta_i).$$ We denote by $M$ the
manifold $B^N$, endowed with a Riemannian metric $ds^2$ such that
on $\Phi_i(V_i)$, for $i>1$,
\begin{equation}\label{metric} ds^2:=f(t) dt^2 + g(t)
d\theta^2,\end{equation}
where $f(t)>0$, $g(t)>0$ for every $t\in (0,+\infty)$ and $d \theta^2$ is
the standard metric on $\Sphere^{N-1}$. $ds^2$ is well-defined on
$B^N \setminus \left\{ 0\right\}$.  \par
We suppose that the metric is asymptotically hyperbolic, that is, as
$t\rightarrow +\infty$,
\begin{equation}\label{condinf} f(t) \rightarrow 1, \quad \quad
g(t)\rightarrow \sinh^2 t.\end{equation} As for the behaviour as
$t\rightarrow 0$, we suppose that for $t\in (0,\epsilon)$
($\epsilon= 2\,{\rm settanh} (\delta)$)
\begin{equation}\label{condzero} f(t)\equiv 1, \quad \quad g(t)=t^2.
\end{equation}
This assures that $ds^2$ can be extended to a smooth Riemannian
metric on all $M$; indeed, for $t\in (0,\epsilon)$, $ds^2$ is the
expression, in polar coordinates, of the Euclidean metric on
$\Real^N$. As already remarked in the Introduction, the essential
spectrum of the Laplace-Beltrami operator acting on $p$-forms on a
complete noncompact Riemannian manifold does not change under
perturbations of the Riemannian metric on compact sets
(\cite{Eichhorn}). As a consequence, condition (\ref{condzero})
does not modify essentially the spectral properties of the
Laplace-Beltrami operator on $M$.\par The manifold $M$, endowed
with the Riemannian metric  $ds^2$, is complete. Indeed, in view
of (\ref{condzero}) and (\ref{condinf}), there exist $C_1,C_2>0$,
$D_1,D_2>0$ such that for every $t>0$ $$C_1\leq f(t) \leq C_2, $$
$$ D_1 \sinh^2 t \leq g(t) \leq D_2 \sinh^2 t;$$ hence the
distance $d_M$ induced by $ds^2$, given by
$$d_M(p_1,p_2)=\inf_{\gamma \in \Gamma(p_1,p_2)}\int_0^1
\left(f(t(s))\left(\frac{d\gamma^1}{ds}\right)^2+g(t(s))
\|\frac{d\gamma^i}{ds}\|^2_{\Sphere^{N-1}}\right)^{\frac{1}{2}}
\,ds,$$ is equivalent to the distance induced by the hyperbolic
metric, which is complete. \par For $p=0,...,N$, we will denote by
$C^{\infty}(\Lambda^p(M))$ the space of all smooth $p$-forms on
$M$, and by $C^{\infty}_c(\Lambda^p(M))$ the set of all smooth,
compactly supported $p$-forms on $M$. For any $\omega\in
C^{\infty}(\Lambda^p(M))$, we will denote by $|\omega(t,\theta)|$
the norm induced by the Riemannian metric on the fiber over
$(t,\theta)$, given in local coordinates by
$$|\omega(t,\theta)|^2=g^{i_1j_1}(t,\theta)...g^{i_pj_p}(t,\theta)
\omega_{i_1...i_p}(t,\theta)\omega_{j_1...j_p}(t, \theta),$$ where
$g^{ij}$ is the expression of the Riemannian metric in local
coordinates. We will denote by $d_M$, $*_M$, $\delta_M$,
respectively, the differential, the Hodge $*$ operator and the
codifferential on $M$, defined as in \cite{deRham}. $\Delta_M$ will stand
for
the
Laplace-Beltrami operator acting on $p$-forms $$
\Delta_M=d_M\delta_M+\delta_M d_M,$$ which is expressed in local
coordinates by the Weitzenb\"ock formula
$$((\Delta_M)\omega)_{i_1...i_p}= - g^{ij}\nabla_i \nabla_j
\omega_{i_1...i_p}+\sum_j R^{\alpha}_{j}
\omega_{i_1...\alpha...i_p} + \sum_{j,l\not=j}
R^{\alpha\,\beta\,}_{\,i_j\,i_l} \omega_{\alpha
i_1...\beta...i_p}, $$ where $\nabla_i\omega$ is the covariant
derivative of $\omega$ with respect to the Riemannian metric, and
$R^i_j$, $R^{i\,j\,}_{\,k\,l}$ denote respectively the local
components of the Ricci tensor and the Riemann tensor induced by
the Riemannian metric. As usual, $L^2_p(M)$ will denote the
completion of $C^{\infty}_c(\Lambda^p(M))$ with respect to the
norm $\|\omega\|_{L^2_p(M)}$ induced by the scalar product
$$\langle\omega,\tilde\omega\rangle_{L^2_p(M)}:= \int_M \omega
\wedge
*_M\tilde \omega
; $$
$\|\omega\|_{L^2_p(M)}$ reads also
$$ \|\omega\|^2_{L^2_p(M)}= \int_M |\omega(t,\theta)|^2 dV_M,$$
where $dV_M$ is the volume element of $(M,ds^2)$.\par
It is well-known that, since the Riemannian metric on $M$ is complete, the
Laplace-Beltrami operator is essentially selfadjoint on
$C^{\infty}_c(\Lambda^p(M))$, for $p=0,...,N$. We will denote by
$\Delta_M$ also its closure.\par
Now, given $\omega \in C^{\infty}(\Lambda^p(M))$, let us write
\begin{equation}\label{decom1} \omega= \omega_1 + \omega_2
\wedge dt,\end{equation}
where $\omega_1$ and $\omega_2$ are respectively a $p$-form and a
$(p-1)$-form on $\Sphere^{N-1}$ depending on $t$. An easy
computation shows that $*_M \omega$ can be expressed in terms of
(\ref{decom1}) as
\begin{multline}\label{Hodge1}*_M \omega=
(-1)^{N-p}g^{\frac{N-2p+1}{2}}(t) f^{-\frac{1}{2}}(t)
*_{\Sphere^{N-1}} \omega_2 +\\
g^{\frac{N-2p-1}{2}}(t)f^{\frac{1}{2}}(t)*_{\Sphere^{N-1}}\omega_1\wedge
dt, \end{multline}
where $*_{\Sphere^{N-1}}$ denotes the Hodge $*$ operator on
$\Sphere^{N-1}$. Moreover, $d_M$ and $\delta_M$ split respectively as
\begin{equation}\label{dM1}d_M\omega= d_{\Sphere^{N-1}}\omega_1 + \left\{
(-1)^p \frac{\partial \omega_1}{\partial t}+ d_{\Sphere^{N-1}}\omega_2
\right\}\wedge dt ,\end{equation}
\begin{multline}\label{deltaM1}\delta_M \omega= g^{-1}(t)
\delta_{\Sphere^{N-1}} \omega_1 + (-1)^p
f^{-\frac{1}{2}}
g^{\frac{-N-1+2p}{2}}\frac{\partial}{\partial t}
\left(f^{-\frac{1}{2}} g^{\frac{N+1-2p}{2}}\omega_2\right)+\\+ g^{-1}
\delta_{\Sphere^{N-1}}\omega_2 \wedge dt,\end{multline}
where $p$ is the degree of $\omega$, $d_{\Sphere^{N-1}}$ is the
differential on $\Sphere^{N-1}$ and $\delta_{\Sphere^{N-1}}$ is the
codifferential on $\Sphere^{N-1}$. \par
Moreover, the $L^2$-norm of $\omega \in
C^{\infty}(\Lambda^p(M))\cap L^2_p(M)$ can be written as
\begin{multline}\label{norm1} \|\omega\|^2_{L^2_p(M)}= \int_0^{+\infty}
g^{\frac{N-2p-1}{2}}(s) f^{\frac{1}{2}}(s)
\|\omega_1(s)\|^2_{L^2_p(\Sphere^{N-1})} \,ds +\\
+ \int_0^{+\infty} g^{\frac{N+1-2p}{2}}(s)f^{-\frac{1}{2}}(s)
\|\omega_2(s)\|^2_{L^2_{p-1}(\Sphere^{N-1})} \,ds,
\end{multline}
where $\|.\|_{L^2_p(\Sphere^{N-1})}$ is the $L^2$-norm for $p$-forms on
$\Sphere^{N-1}$.
\section{Zero in the spectrum}
In the present section we will investigate whether $0$ belongs or
not to the point (and essential) spectrum of $\Delta_M$, for
differential forms of degree $p=0,...,N$. The main tool employed
is the following generalization of a result of Dodziuk
(\cite{Dodziuk}):
\begin{thm}
Let us consider, for $N\geq 2$, the manifold $M$ endowed
with a complete Riemannian metric of type (\ref{metric}),
satisfying condition (\ref{condzero}) for $t\in (0,\epsilon)$ ; then, if
we denote by ${\mathcal H}^p(M)$, for $p=0,...,N$,
the space of $L^2$ harmonic $p$-forms on $M$, we have
\begin{enumerate}
\item for $p\notin \left\{0,N,N/2\right\}$, ${\mathcal
H}^p(M)=\left\{0\right\}$;
\item if $\int_0^{\infty}f^{\frac{1}{2}}(s)g^{\frac{N-1}{2}}(s)\,ds=
+\infty$,
${\mathcal H}^N(M)\simeq {\mathcal H}^0(M)=\left\{0\right\}$; if on the
contrary $\int_0^{\infty}f^{\frac{1}{2}}(s)g^{\frac{N-1}{2}}(s)\,ds<+
\infty$,
${\mathcal H}^N(M)\simeq {\mathcal H}^0(M)=\Real$;
\item if $p=\frac{N}{2}$, ${\mathcal H}^p(M)= \left\{0\right\}$ if
$\int_1^{+\infty}f^{\frac{1}{2}}(s)g^{-\frac{1}{2}}(s)\,ds = +\infty$; if
on the other hand
$\int_1^{+\infty}f^{\frac{1}{2}}(s)g^{-\frac{1}{2}}(s)\,ds < +\infty$,
${\mathcal H}^{\frac{N}{2}}(M)$ is a Hilbert space of infinite dimension.
\end{enumerate}
\end{thm}
\begin{proof}
The proof follows very closely the argument in \cite{Dodziuk}; it will be
exposed here for the sake of completeness.\par
An $L^2$-form
on $M$ is
harmonic if and only it is closed and coclosed. Hence, $\omega \in
{\mathcal H}^p(M)$ if and only if
\begin{equation}\label{harm} \|\omega\|_{L^2_p(M)}<\infty, \quad \quad
d\omega=0,
\quad \quad d *_M \omega=0. \end{equation}
Moreover, $*_M$ gives an isomorphism between
${\mathcal H}^p(M)$ and ${\mathcal H}^{N-p}(M)$.\par
The proof of 2) is immediate; if $\omega$ is a harmonic function, not
identically vanishing,
$\omega$ is constant on $M$, hence $\omega \in L^2(M)$ if and only if
the
total volume of $M$, given by
$\int_0^{\infty}f^{\frac{1}{2}}(s)g^{\frac{N-1}{2}}(s)\,ds$, is finite.
\par
We now come to the proof of 1). Let $\omega\in {\mathcal H}^p(M)$,
for $p\not=  0,N$, and let us consider its decomposition
(\ref{decom1}). Then, in view of (\ref{dM1}), $d_M \omega=0$
implies $$d_{\Sphere^{N-1}}\omega_1=0,\quad \quad
d_{\Sphere^{N-1}}\omega_2 + (-1)^p \frac{\partial
\omega_1}{\partial t}=0, $$ whilst $d_M*_M \omega=0$ yields
$$d_{\Sphere^{N-1}}*_{\Sphere^{N-1}}\omega_2=0, $$
\begin{equation}\label{zero}
g^{\frac{N-2p-1}{2}}f^{\frac{1}{2}}d_{\Sphere^{N-1}}*_{\Sphere^{N-1}}
\omega_1 + \frac{\partial}{\partial t} \left(
g^{\frac{N-2p+1}{2}}f^{\frac{1}{2}}*_{\Sphere^{N-1}}
\omega_2\right)=0.
\end{equation}
In view of (\ref{norm1}), the boundedness of the $L^2$-norm of $\omega$
reads
\begin{multline} \int_0^{+\infty}\int_{\Sphere^{N-1}}
(g^{\frac{N-2p-1}{2}}
f^{\frac{1}{2}}|\omega_1(t,\theta)|^2\\+g^{\frac{N-2p+1}{2}}
f^{-\frac{1}{2}}|\omega_2(t,\theta)|^2 )\,dV_{\Sphere^{N-1}}\,dt <
+\infty;\end{multline} moreover, since $|\omega(t,\theta)|$ is
bounded in a neighbourhood of $0$, we have that $$
|\omega(t,\theta)|^2 = g(t)^{-p}|\omega_1(t,\theta)|^2 +
f(t)^{-1}g(t)^{1-p}|\omega_2(t,\theta)|^2 \leq C$$ for some $C>0$
for $t\in (0,\epsilon]$. \par Applying $*_{\Sphere^{N-1}}$ to both
sides of (\ref{zero}), we find the following set of conditions:
\begin{equation}\label{zero1} d_{\Sphere^{N-1}} \omega_1 =0;
\end{equation}
\begin{equation}\label{zero2} d_{\Sphere^{N-1}}
*_{\Sphere^{N-1}}\omega_2=0;\end{equation}
\begin{equation}\label{zero3} d_{\Sphere^{N-1}}\omega_2 +(-1)^p
\frac{\partial \omega_1}{\partial t}=0; \end{equation}
\begin{equation}\label{zero4} \frac{\partial}{\partial
t}\left(g^{\frac{N-2p+1}{2}}(t)f^{-\frac{1}{2}}(t) \omega_2\right)+
(-1)^p f^{\frac{1}{2}}(t) g^{\frac{N-2p-1}{2}}(t)
\delta_{\Sphere^{N-1}}\omega_1=0; \end{equation}
\begin{equation}\label{zero5} g^{-p}(t) |\omega_1 (t,\theta)|^2 +
f^{-1}(t) g^{1-p}(t)|\omega_2(t,\theta)|^2 \leq C \quad \forall t\in
(0,\epsilon];
\end{equation}
\begin{multline}\label{zero6} \int_0^{+\infty}\int_{\Sphere^{N-1}}
(g^{\frac{N-2p-1}{2}} f^{\frac{1}{2}}|\omega_1(t,\theta)|^2+\\ +
g^{\frac{N-2p+1}{2}}f^{-\frac{1}{2}}|\omega_2(t,\theta)|^2
)\,dV_{\Sphere^{N-1}}\,dt < +\infty.\end{multline} Now, it can be
shown that if $\omega \in {\mathcal H}^p(M)$ and $\omega_1=0$,
then $\omega_2=0$; indeed, if $\omega_2\wedge dt \in {\mathcal
H}^p(M)$, in view of (\ref{zero2}) and (\ref{zero3}) $\omega_2$ is
a harmonic form on $\Sphere^{N-1}$ for every $t>0$. Since $0 \leq
p-1 \leq N-2$, $\omega_2(t,\theta)$ can be nonzero only if $p-1=
\deg \omega_2=0$, that is, only if $\omega_2$ is a function not
depending on $\theta$. On the other hand, (\ref{zero4}) implies $$
\frac{\partial}{\partial t}\left(g^{\frac{N-1}{2}}f^{-\frac{1}{2}}
\omega_2\right)=0,$$ that is, $\omega_2= C
g(t)^{-\frac{N-1}{2}}f(t)^{\frac{1}{2}}$, which diverges as $t
\rightarrow 0$, in contradiction with (\ref{zero5}), unless
$C=0$.\par Hence, if $\omega\not=0$ and $\omega \in {\mathcal
H}^p(M)$, then $\omega_1 \not=0$. Now, applying
$d_{\Sphere^{N-1}}$ to both sides of (\ref{zero4}), since
$d_{\Sphere^{N-1}}$ commutes with $ \frac{\partial}{\partial t}$,
we get $$\frac{\partial}{\partial t}
\left(g(t)^{\frac{N-2p+1}{2}}f^{-\frac{1}{2}} d_{\Sphere^{N-1}}
\omega_2 \right)+ (-1)^p
f(t)^{\frac{1}{2}}g(t)^{\frac{N-2p-1}{2}}d_{\Sphere^{N-1}}
\delta_{\Sphere^{N-1}} \omega_1=0, $$ whence, in view of
(\ref{zero3}), $$\frac{\partial}{\partial
t}\left(g(t)^{\frac{N-2p+1}{2}}f(t)^{-\frac{1}{2}}\frac{\partial\omega_1}
{\partial t} \right)=
f(t)^{\frac{1}{2}}g(t)^{\frac{N-2p-1}{2}}d_{\Sphere^{N-1}}
\delta_{\Sphere^{N-1}} \omega_1. $$ Taking, for fixed $t>0$, the
scalar product of both sides of the last equation with $\omega_1$,
we get $$\langle\frac{\partial}{\partial
t}\left(g(t)^{\frac{N-2p+1}{2}}f(t)^{-\frac{1}{2}}\frac{\partial\omega_1}
{\partial t} \right), \omega_1\rangle_{L^2_p(\Sphere^{N-1})}$$ $$=
\langle\delta_{\Sphere^{N-1}}\omega_1,\delta_{\Sphere^{N-1}}
\omega_1\rangle_{L^2_p(\Sphere^{N-1})} \geq 0,$$ whence
$$\frac{\partial}{\partial t} \langle
g(t)^{\frac{N-2p+1}{2}}f(t)^{\frac{1}{2}}\frac{\partial
\omega_1}{\partial t},\omega_1\rangle_{L^2_p(\Sphere^{N-1})}=$$
$$= \langle\frac{\partial}{\partial
t}\left(g(t)^{\frac{N-2p+1}{2}}f(t)^{-\frac{1}{2}} \frac{\partial
\omega_1}{\partial
t}\right),\omega_1\rangle_{L^2_p(\Sphere^{N-1})} $$ $$ +
g(t)^{\frac{N-2p+1}{2}}f(t)^{-\frac{1}{2}}\langle \frac{\partial
\omega_1}{\partial t}, \frac{\partial \omega_1}{\partial
t}\rangle_{L^2_p(\Sphere^{N-1})} \geq 0.$$ Due to the boundedness
of $|\omega|$ near $0$ and to (\ref{condzero}),
$|\omega_1(t,\theta)|_{\Sphere^{N-1}}=O(t^{2p})$ for small $t$. As
a consequence, $$\langle
f(t)^{-\frac{1}{2}}g(t)^{\frac{N-2p+1}{2}}\frac{\partial
\omega_1}{\partial t},\omega_1\rangle_{L^2_p(\Sphere^{N-1})}=
O(t^N),$$ hence $$ \frac{\partial}{\partial t}
\langle\omega_1,\omega_1\rangle_{L^2_p(\Sphere^{N-1})}= 2\,
\langle\frac{\partial \omega_1}{\partial
t},\omega_1\rangle_{L^2_p(\Sphere^{N-1})}\geq 0$$ for every $t>0$,
that is, $\|\omega_1(t)\|_{L^2_p(\Sphere^{N-1})}$ is a
nondecreasing function of $t$.\par Now, let $\omega_1 \not=0$;
since $\|\omega_1(t)\|_{L^2_p(\Sphere^{N-1})}$ is nondecreasing
and \\ $\|\omega\|_{L^2_p(M)} < +\infty$, $$ \int_1^{+\infty}
g(s)^{\frac{N-2p-1}{2}}f(s)^{\frac{1}{2}}\,ds$$ $$ \leq C
\int_1^{+\infty} g(s)^{\frac{N-2p-1}{2}}f(s)^{\frac{1}{2}}
\|\omega_1(s)\|^2_{L^2_p(\Sphere^{N-1})}\,ds \leq
\|\omega\|^2_{L^2_p(M)}<+\infty.$$ Hence for $p\not=0,N$,
${\mathcal H}^p(M) \not= \left\{0\right\}$ implies
$$\int_1^{+\infty} g(s)^{\frac{N-2p-1}{2}}f(s)^{\frac{1}{2}}\,ds<
+\infty, $$ and, by duality, $$\int_1^{+\infty}
g(s)^{\frac{-N+2p-1}{2}}f(s)^{\frac{1}{2}}\,ds < +\infty.$$ If
$N=2p$, the two integrands coincide. If, on the contrary, $N-2p
\not=0$, then, since $ (N-2p-1)(-N+2p-1)= 1-(N-2p)^2$, either one
of the exponents is zero, or the two exponents have opposite
signs; in both cases one of the integrals diverges. Hence, for
$p\notin \left\{ 0,N,N/2\right\}$, ${\mathcal H}^p(M) = \left\{
0\right\}$.\par Finally we come to 3). For $p= N/2$, if
$\int_1^{+\infty}g(s)^{-1/2}f(s)^{1/2}\,ds =+\infty,
$
${\mathcal H}^{p}(M)=\left\{0\right\}$. This proves the first half of 3).
We still have to prove that if $\int_1^{+\infty}g(s)^{-1/2}f(s)^{1/2}\,ds
<+\infty$, ${\mathcal H}^{N/2}(M)$ has infinite dimension.
To this purpose, let us recall that if $N=2p$ the Hodge $*$
operator  acting on forms of degree $p$ depends only on the
conformal structure of the manifold. Hence the conditions
$\|\omega\|_{L^2_p}<+\infty$, $d \omega =0$, $d*\omega=0$ are
conformally invariant.\par Now, let us suppose that
$\int_1^{+\infty} g(s)^{-1/2}f(s)^{1/2}\,ds
<+\infty$, and let us denote by $B(0,r)$ the open ball in $\Real^N$ with
radius
$$r=\exp\left(\int_1^{+\infty} g(s)^{-1/2}f(s)^{1/2}\,ds\right)$$ centered
in
$0$,
endowed with polar coordinates. Then
consider the mapping:
$$F:M\setminus \left\{0\right\} \longrightarrow
\Real^N\setminus
\left\{0\right\} $$
given by
$$F(t,\theta):= \left(\exp \left(\int_1^t
g(s)^{-1/2}f(s)^{1/2}\,ds\right),
\theta\right). $$
In view of condition (\ref{condzero}), $F$ can be extended to a
$C^1$-diffeomorphism of $M$ into $B(0,r)$, which is
actually $C^{\infty}$ on $M\setminus \left\{0\right\}$.
Moreover,
an easy computation shows that $F$ is conformal from $M$, endowed with the
metric (\ref{metric}), to $B(0,r)$, endowed
with the Euclidean metric.\par
Let us denote by ${\mathcal H}$
the (infinite-dimensional) space of all smooth $p$-forms on
$B(0,r)$ harmonic with
respect to the Euclidean metric;
since $F$ is conformal and $N=2p$, $F^* {\mathcal H}$ consists of forms of
degree $p$, square-summable on $M$, smooth on $M$ (up to modifications at
$0$) and harmonic. As a consequence, ${\mathcal H}^{N/2}(M)$ has infinite
dimension.
\end{proof}
In our case, since $ f(t)\rightarrow 1$ and $g(t) \rightarrow
\sinh^2 t$ as $t\rightarrow +\infty$, then
$$\int_0^{+\infty}f(s)^{\frac{1}{2}}g(s)^{\frac{N-1}{2}}\,ds
=+\infty,$$ whilst
$$\int_1^{+\infty}f(s)^{\frac{1}{2}}g(s)^{-\frac{1}{2}}\,ds<
+\infty.$$ As a consequence we can easily deduce the following
\begin{thm}\label{zeroess} For $N\geq 2$,  let us consider the manifold
$M$, endowed with
a Riemannian metric of type (\ref{metric}), satisfying conditions
(\ref{condinf}) and (\ref{condzero}). Then \begin{enumerate}
\item if $p\not=N/2$, then $0\notin
\sigma_p(\Delta_M)$;
\item if $p=N/2$, ${\mathcal H}^p(M)$ is a Hilbert space of infinite
dimension, hence $0\in \sigma_{{\rm ess}}(\Delta_M)\cap
\sigma_p(\Delta_M)$.
\end{enumerate}
\end{thm}

\section{Hodge decomposition and unitary equivalence}
From (\ref{dM1}) and (\ref{deltaM1}), a
lengthy but straightforward computation gives
$$ \Delta_M \omega= (\Delta_M \omega)_1 + (\Delta_M \omega)_2\wedge dt,$$
where
\begin{multline}\label{LB11}(\Delta_M \omega)_1=
g^{-1}(t)\Delta_{\Sphere^{N-1}} \omega_1 + (-1)^p f^{-1}(t)
g^{-1}(t) \frac{\partial g}{\partial t} d_{\Sphere^{N-1}}\omega_2 +\\
-f^{-\frac{1}{2}}(t) g^{\frac{-N+1+2p}{2}}(t)
\frac{\partial}{\partial t} \left( f^{-\frac{1}{2}}(t)
g^{\frac{N-1-2p}{2}}(t) \frac{\partial \omega_1}{\partial
t} \right)
\end{multline} and
\begin{multline}\label{LB12}(\Delta_M \omega)_2 = g^{-1}(t)
\Delta_{\Sphere^{N-1}} \omega_2  + (-1)^{p}g^{-2}(t)
\frac{\partial g}{\partial t} \delta_{\Sphere^{N-1}}\omega_1  +\\
- \frac{\partial}{\partial t} \left\{f^{-\frac{1}{2}}(t)
g^{\frac{-N-1+2p}{2}}(t) \frac{\partial}{\partial
t}\left(f^{-\frac{1}{2}}(t) g^{\frac{N+1-2p}{2}}(t)
\omega_2\right) \right\}. \end{multline} Here we denote by
$\Delta_{\Sphere^{N-1}}$ the Laplace-Beltrami operator on
$\Sphere^{N-1}$.\par Since for every $\omega\in
C^{\infty}(\Lambda^p(M))\cap L^2_p(M)$ we have that $\omega_1 \in
L^2_p(M)$, $\omega_2 \wedge dt \in L^2_p(M)$ and $$
\langle\omega_1,\omega_2\wedge dt\rangle_{L^2_p(M)}=0,$$
 (\ref{decom1}) gives rise to an orthogonal decomposition of
$L^2_p(M)$ into two closed subspaces. However, (\ref{LB11}) and
(\ref{LB12}) show that $\Delta_M$ is not invariant under this
decomposition. As a consequence, further decompositions are
required.\par It is well-known that, for $0\leq p \leq N-1$, $$
C^{\infty}(\Lambda^p(\Sphere^{N-1}))=
dC^{\infty}(\Lambda^{p-1}(\Sphere^{N-1}))\oplus \delta
C^{\infty}(\Lambda^{p+1}(\Sphere^{N-1})) \oplus {\mathcal
H}^p(\Sphere^{N-1}),$$
where ${\mathcal H}^p(\Sphere^{N-1})$ is the space of harmonic $p$-forms
on $\Sphere^{N-1}$ (empty if $p \not= 0, N-1$), and the decomposition is
orthogonal in $L^2_p(M)$. Hence, for $0\leq p \leq N-1$,
$$L^2_p(\Sphere^{N-1})=
\overline{dC^{\infty}(\Lambda^{p-1}(\Sphere^{N-1}))}\oplus
\overline{\delta C^{\infty}(\Lambda^{p+1}(\Sphere^{N-1}))}\oplus
{\mathcal H}^p(\Sphere^{N-1}). $$
Thus, for $1\leq p \leq N-1$, every $\omega\in L^2_p(M)$ can be written as
\begin{equation}\label{decom2}\omega=\omega_{1\delta}\oplus
\omega_{2d}\wedge dt \oplus
(\omega_{1d} \oplus \omega_{2\delta}\wedge dt),\end{equation} where
$\omega_{1\delta}$ (resp. $\omega_{1d}$) is a coclosed (resp. closed)
$p$-form on $\Sphere^{N-1}$ parametrized by $t$, and $\omega_{2\delta}$
(resp. $\omega_{2d}$) is a coclosed (resp. closed) $(p-1)$-form on
$\Sphere^{N-1}$ parametrized by $t$. In this way we get the orthogonal
decomposition
$$ L^2_p(M)= {\mathcal
L}^1(M)\oplus {\mathcal L}^2(M)\oplus{\mathcal L}^3(M),$$
where for every $\omega\in L^2_p(M)$, $\omega_{1\delta}\in {\mathcal
L}^1(M)$, $\omega_{2d}\wedge dt\in {\mathcal L}^2(M)$ and
$\omega_{1d}\oplus (\omega_{2\delta}\wedge dt)\in {\mathcal L}^3(M)$.
Since
$$d_{\Sphere^{N-1}}\Delta_{\Sphere^{N-1}}=\Delta_{\Sphere^{N-1}}
d_{\Sphere^{N-1}},\quad \quad
\delta_{\Sphere^{N-1}}\Delta_{\Sphere^{N-1}}=\Delta_{\Sphere^{N-1}}
\delta_{\Sphere^{N-1}},$$ $$\frac{\partial}{\partial
t}d_{\Sphere^{N-1}}=d_{\Sphere^{N-1}} \frac{\partial}{\partial t}, \quad
\quad \frac{\partial}{\partial t}\delta_{\Sphere^{N-1}}=
\delta_{\Sphere^{N-1}}
\frac{\partial}{\partial t},$$ the Laplace-Beltrami operator is invariant
under this decomposition, and can be written as the orthogonal sum
$$ \Delta_M= \Delta_{M1}\oplus \Delta_{M2} \oplus \Delta_{M3}.$$
It is easy to see that, for $i=1,2,3$, $\Delta_{Mi}$ is essentially
selfadjoint on $C^{\infty}_c(\Lambda^p(M))\cap {\mathcal L}^i$. We denote
again by $\Delta_{Mi}$ its closure.\par
Since the orthogonal sum is finite, for $1\leq p \leq N-1$,
$$ \sigma_{\rm ess}(\Delta_M)= \bigcup_{i=1}^3 \sigma_{\rm
ess}(\Delta_{Mi}),$$
$$ \sigma_p(\Delta_M)=\bigcup_{i=1}^3 \sigma_p (\Delta_{Mi}).$$
For $p=0$ (resp. $p=N$), any $\omega \in L^2(M)$ can be written as
$\omega = \omega_{1\delta}$ (resp. $\omega= \omega_{2d} \wedge dt$), where
$\omega_{1\delta}$ (resp. $\omega_{2d}$) is a coclosed (resp. closed)
$0$-form (resp. $(N-1)$-form) parametrized by $t$ on $\Sphere^{N-1}$.
Hence $L^2_0(M)= {\mathcal L}^1$ (resp. $L^2_{N-1}(M)= {\mathcal L}^2$)
and $\Delta_M =\Delta_{M1}$ (resp. $\Delta_M = \Delta_{M2}$). \par
As a consequence, in order to determine the spectrum of $\Delta_M$ it
suffices to study the spectral properties of $\Delta_{Mi}$, $i=1,2,3$.
\par
Then, let us introduce a further decomposition. First of all, we
decompose $\omega_{1\delta}$ according to an orthonormal basis
$\left\{\tau_{1k}\right\}_{k\in \mathbb N}$ of coclosed
$p$-eigenforms of $\Delta_{\Sphere^{N-1}}$;
 this yields
\begin{equation}\label{1delta} \omega_{1\delta}=\oplus_{k}
h_{k}(t) \tau_{1k},\end{equation} where $h_{k}(t)\tau_{1k}\in
L^2_p(M)$ for every $k\in \mathbb N$, and the sum is orthogonal in
$L^2_p(M)$, thanks to (\ref{metric}). We will call $p$-form of
type I any $p$-form $\omega\in L^2_p(M)$ such that $$\omega= h(t)
\tau_1,$$ where $\tau_1$ is a coclosed normalized $p$-eigenform of
$\Delta_{\Sphere^{N-1}}$, corresponding to some eigenvalue
$\lambda$. For every $k\in \mathbb N$, let us denote by
$\lambda_k^p\in \sigma_p(\Delta_{\Sphere^{N-1}})$ the eigenvalue
associated to $\tau_{1k}$. Since for every $k\in \mathbb N$
\begin{multline}\label{h1}\Delta_{M1}(h(t)\tau_{1k})=
\frac{\lambda_k^p}{g(t)} h(t)\tau_{1k} \\-
f(t)^{-\frac{1}{2}}g(t)^{\frac{-N+1+2p}{2}}\frac{\partial}{\partial
t}\left(f(t)^{-\frac{1}{2}}g(t)^{\frac{N-1-2p}{2}}\frac{\partial
h}{\partial
t}\right)
\tau_{1k}, \end{multline}
$\Delta_{M1}$ is invariant under the decomposition (\ref{1delta}), and,
since
if $\omega=h(t)\tau_{1k}$
$$\|\omega\|^2_{L^2_p(M)}=
\int_0^{\infty}g(s)^{\frac{N-2p-1}{2}}f(s)^{\frac{1}{2}}h(s)^2\,ds,$$
$\Delta_{M1}$ is unitarily equivalent to the direct sum with respect to
$k\in \mathbb N$
of
the operators
$$\Delta_{1\lambda_k^p}:{\mathcal D}(\Delta_{1\lambda_k^p})\subset
L^2(\Real^+, g^{\frac{N-2p-1}{2}}f^{\frac{1}{2}})\longrightarrow
L^2(\Real^+, g^{\frac{N-2p-1}{2}}f^{\frac{1}{2}})$$
\begin{multline} \Delta_{1\lambda_k^p}h=\left\{
\frac{\lambda_k^p}{g(t)} h(t)  -
f(t)^{-\frac{1}{2}}g(t)^{\frac{-N+1+2p}{2}}\frac{\partial}{\partial
t}\left(f(t)^{-\frac{1}{2}}g(t)^{\frac{N-1-2p}{2}}\right)
\right\}.\end{multline} If we introduce the transformation
\begin{equation}\label{trasf1}w(t)=
h(t)f(t)^{\frac{1}{4}}g(t)^{\frac{N-2p-1}{4}},
\end{equation}
a direct (but lengthy) computation shows that $\Delta_{M1}$ is
unitarily equivalent to the direct sum, over $k\in \mathbb N$, of
the operators $$D_{1\lambda_k^p}: {\mathcal
D}(D_{1\lambda_k^p})\subset L^2(\Real^+)\longrightarrow
L^2(\Real^+)$$ given by
\begin{multline}\label{w1} D_{1\lambda_k^p}w = - \frac{\partial}{\partial
t}\left( \frac{1}{f} \frac{\partial w}{\partial t} \right)+
\left\{ -\frac{7}{16} \frac{1}{f^3} \left(\frac{\partial
f}{\partial t}\right)^2+ \frac{1}{4} \frac{1}{f^2}
\frac{\partial^2 f}{\partial t^2}\right.\\ - \frac{1}{2}
\frac{1}{f^2}\frac{\partial f}{\partial t}\frac{(N-1-2p)}{4}
\frac{1}{g}\frac{\partial g}{\partial t} + \frac{1}{f}
\frac{(N-2p-1)}{4}
\frac{(N-2p-5)}{4}\frac{1}{g^2}\left(\frac{\partial g}{\partial
t}\right)^2\\ \left. +
\frac{1}{f}\frac{(N-2p-1)}{4}\frac{1}{g}\frac{\partial^2 g
}{\partial t^2} + \frac{\lambda_k^p}{g}\right\} w.
\end{multline}
Analogously, we decompose $\omega_{2d}$ according to an
orthonormal basis of closed $(p-1)$-eigenforms $\left\{\tau_{2k}
\right\}_{k\in \mathbb N}$ of $\Delta_{\Sphere^{N-1}}$:
\begin{equation}\label{2d} \omega_{2d}\wedge dt=\oplus_{k}
h_{k}(t) \tau_{2k}\wedge dt.\end{equation} We will call $p$-form
of type II a $p$-form $\omega\in L^2_p(M)$ such that $$\omega=
h(t) \tau_2 \wedge dt,$$ where $\tau_2$ is a coclosed normalized
$(p-1)$-eigenform, corresponding to some eigenvalue $\lambda$ of
$\Delta_{\Sphere^{N-1}}$. For every $k\in \mathbb N$
$$\Delta_{M2}(h(t)\tau_{2k}\wedge dt )= (\Delta_{2\lambda_k^{p-1}}
h )\tau_{2k}\wedge dt,$$ where
\begin{multline}\label{h2}
\Delta_{2\lambda_k^{p-1}} h = \frac{\lambda_k^{p-1}}{g(t)} h(t)-
\\ \frac{\partial}{\partial t}\left\{
f(t)^{-\frac{1}{2}}g(t)^{\frac{-N-1+2p}{2}}
\frac{\partial}{\partial t}\left( f(t)^{-\frac{1}{2}}
g(t)^{\frac{N+1-2p}{2}}h(t) \right) \right\}. \end{multline} Here,
again, for every $k\in \mathbb N$ we denote by $\lambda_k^{p-1}$
the eigenvalue of $\Delta_{\Sphere^{N-1}}$ corresponding to the
eigenform $\tau_{2k}$. Since if $\omega= h(t) \tau_{2k}\wedge dt$
$$\|\omega\|^2_{L^2_p(M)}=
\int_0^{\infty}g(s)^{\frac{N-2p+1}{2}}f(s)^{-\frac{1}{2}}h(s)^2\,ds,$$
introducing the transformation
\begin{equation}\label{trasf2} w(t)= h(t)
f(t)^{-\frac{1}{4}}g(t)^{\frac{N+1-2p}{4}},
\end{equation}
we find that $\Delta_{2M}$ is unitarily equivalent to the direct sum, with
respect to
$k\in \mathbb N$, of the operators
$$D_{2 \lambda^{p-1}_k}:{\mathcal D}(D_{2 \lambda^{p-1}_k})\subset
L^2(\Real^+)
\longrightarrow L^2(\Real^+) $$
\begin{multline}\label{w2} D_{2\lambda^{p-1}_k}w = -
\frac{\partial}{\partial t}\left( \frac{1}{f} \frac{\partial
w}{\partial t} \right)+  \left\{ -\frac{7}{16} \frac{1}{f^3}
\left(\frac{\partial f}{\partial t}\right)^2+\frac{1}{4}
\frac{1}{f^2} \frac{\partial^2 f}{\partial t^2} \right.\\ -
\frac{1}{2} \frac{1}{f^2}\frac{\partial f}{\partial
t}\frac{(N-1+2p)}{4} \frac{1}{g}\frac{\partial g}{\partial t}+
\frac{1}{f} \frac{(N-2p+1)}{4}
\frac{(N-2p+5)}{4}\frac{1}{g^2}\left(\frac{\partial g}{\partial
t}\right)^2+ \\ \left.
+\frac{1}{f}\frac{(-N+2p-1)}{4}\frac{1}{g}\frac{\partial^2
g}{\partial t^2} + \frac{\lambda_k^{p-1}}{g}\right\} w.
\end{multline}
Finally, we decompose $\omega_{2\delta}$ with respect to an
orthonormal basis of coclosed $(p-1)$-eigenforms $\left\{
\tau_{3k}\right\}_{k\in \mathbb N}$ of $\Delta_{\Sphere^{N-1}}$.
For every $k\in \mathbb N$ we denote by $\lambda_k^{p-1}$ the
eigenvalue corresponding to the eigenform $\tau_{3k}$; then
$\left\{\frac{1}{\sqrt{\lambda_k^{p-1}}}
d_{\Sphere^{N-1}}\tau_{3k} \right\}_{k\in \mathbb N}$ is an
orthonormal basis of closed eigenforms for exact $p$-forms on
$\Sphere^{N-1}$. Hence, we get the following decomposition for
$\omega_{1d}+\omega_{2\delta}\wedge dt$:
\begin{multline}\label{3}\omega_{1d}+\omega_{2\delta}\wedge dt =
\oplus_{k}
\left(\frac{1}{\sqrt{\lambda_k^{p-1}}} h_{1k}
d_{\Sphere^{N-1}}\tau_{3k} \oplus  (-1)^p
h_{2k}\tau_{3k}\wedge dt\right).
\end{multline}
We call $p$-form of type III any $p$-form $\omega$ such that
$$\omega= \frac{1}{\sqrt{\lambda}}h_1(t)d_{\Sphere^{N-1}}\tau_3
\oplus_M (-1)^p h_2(t) \tau_3 \wedge dt ,$$
where $\tau_3$ is a normalized coclosed $(p-1)$-eigenform of
$\Delta_{\Sphere^{N-1}}$, corresponding to the eigenvalue $\lambda$.
A direct computation shows that
\begin{multline}\label{h1h2}
\Delta_{M3}\left(\frac{1}{\sqrt{\lambda^{p-1}_k}}
h_1(t)d_{\Sphere^{N-1}}\tau_{3k}
\oplus_M (-1)^p h_2(t) \tau_{3k} \wedge dt\right)=\\
=\left(\Delta_{1\lambda_k^{p-1}} h_1 + \frac{1}{f(t)}\frac{1}{g(t)}
\frac{\partial
g}{\partial t}\sqrt{\lambda_k^{p-1}} h_2
\right) \left( \frac{1}{\sqrt{\lambda_k^{p-1}}} d_{\Sphere^{N-1}}
\tau_{3k} \right)\\
\oplus \left(\Delta_{2\lambda_k^{p-1}} h_2 + \frac{1}{g^2(t)}
\frac{\partial g}{\partial t}
\sqrt{\lambda_k^{p-1}} h_1 \right) \left((-1)^p \tau_{3k} \wedge dt
\right);
\end{multline}
moreover, if $\omega=\frac{1}{\sqrt{\lambda_k^{p-1}}}
h_1(t)d_{\Sphere^{N-1}}\tau_{3k} \oplus_M (-1)^p h_2(t) \tau_{3k}
\wedge dt$, then $$\|\omega\|^2_{L^2_p(M)}= \int_0^{+\infty}
g(s)^{\frac{N-2p-1}{2}}f(s)^{\frac{1}{2}} h_1(s)^2 \,ds + $$ $$+
\int_0^{+\infty} g(s)^{\frac{N+1-2p}{2}}f(s)^{-\frac{1}{2}}
h_2(s)^2\,ds. $$ Hence, introducing the transformation
\begin{equation}\label{trasf12}\begin{array}{ll}
w_1(t)=&g^{\frac{N-2p-1}{4}}(t)f^{\frac{1}{4}}(t)h_1(t) \\
w_2(t)=& g^{\frac{N-2p+1}{4}}(t)f^{-\frac{1}{4}}(t) h_2(t),
\end{array}
\end{equation}
we find that $\Delta_{M3}$ is unitarily equivalent to the direct sum, with
respect to
$k\in \mathbb N$, of the operators
$$D_{3\lambda_k^{p-1}}:{\mathcal D}(D_{3\lambda^{p-1}_k})\subset
L^2(\Real^+)\oplus
L^2(\Real^+)
\longrightarrow
L^2(\Real^+)\oplus L^2(\Real^+) $$
\begin{multline}\label{w1w2}D_{3\lambda^{p-1}_k}(w_1\oplus w_2)=
\left(D_{1\lambda^{p-1}_k}w_1 +
g(t)^{-\frac{3}{2}}f(t)^{-\frac{1}{2}}\frac{\partial g}{\partial t}
\sqrt{\lambda^{p-1}_k}w_2 \right) \oplus \\ \oplus
\left(D_{2\lambda^{p-1}_k}
w_2 + g(t)^{-\frac{3}{2}}f(t)^{-\frac{1}{2}}\frac{\partial g}{\partial t}
\sqrt{\lambda^{p-1}_k}w_1
\right).
\end{multline}
\par \bigskip As remarked in the Introduction,  J. Eichhorn proved
in \cite{Eichhorn} that for a complete Riemannian metric over a
noncompact manifold the essential spectrum of $\Delta_M$ coincides
with the essential spectrum of the Friedrichs extension
$\Delta_M^F$ of the restriction of $\Delta_M$ to any exterior
domain in $M$.  Thus, if we consider, for $0<\eta<1$, the
Friedrichs extension $\Delta_{M,\eta}^F$ of the operator
$$\Delta'_{M,\eta} : C^{\infty}_c(\Lambda^p(M\setminus B(0,\eta)))
\longrightarrow L^2(M\setminus B(0,\eta))$$ $$ \Delta'_{M,\eta}
\omega= \Delta_M \omega,$$ we have that
 $$\sigma_{\rm
ess}(\Delta_M)= \sigma_{\rm ess}({\Delta_{M,\eta}^F})$$ for every
$\eta$, $0<\eta<1$. Hence, in order to compute the essential
spectrum of $\Delta_M$ it suffices to determine the essential
spectrum of $\Delta_{M,\eta}^F$ for some $\eta$, $0<\eta<1$. For
the sake of simplicity we will write $\Delta_M^F$ instead of
$\Delta_{M,\eta}^F$.
\par The same orthogonal decompositions obtained for $\Delta_M$ hold
also for $\Delta_M^F$: namely, we have a decomposition $$
L^2_p(M\setminus B(0,\eta)) = {\mathcal L}^1(M\setminus
B(0,\eta))\oplus{\mathcal L}^2(M\setminus B(0,\eta)) \oplus
{\mathcal L}^3(M\setminus B(0,\eta))$$ analogous to
(\ref{decom2}), and $\Delta_M^F$ splits accordingly as $$
\Delta_M^F= \Delta_{M1}^F \oplus \Delta_{M2}^F\oplus
\Delta_{M3}^F,$$ where, for $i=1,2,3$, $\Delta_{Mi}^F$ is the
Friedrichs extension of the restriction of $\Delta_M$ to
$C^{\infty}_c(\Lambda^p(M\setminus B(0,\eta))) \cap {\mathcal
H}^i(M\setminus B(0,\eta))$. Moreover (see \cite{Eichhorn}), for
$i=1,2,3$, $\sigma_{\rm ess}(\Delta_{Mi})= \sigma_{\rm
ess}(\Delta_{Mi}^F)$.\par Let $c= {\rm settanh} (\eta)$; again, it
is possible to show that, for $i=1,2$, $\Delta_{Mi}^F$ is
unitarily equivalent to the direct sum, over $k\in \mathbb N$, of
the Friedrichs extensions $D^F_{i\lambda_k^p}$ of the operators
$$D'_{i\lambda_k^p}: C^{\infty}_c(c,+\infty)\longrightarrow
L^2(c,+\infty)$$ given by (\ref{w1}) if $i=1$ and by (\ref{w2}) if
$i=2$.
\par Analogously, $\Delta_{M3}^F$ is unitarily equivalent to the
direct sum, over $k\in \mathbb N$, of the Friedrichs extensions
$D^F_{3\lambda_k^{p-1}}$ of the operators $$
D'_{3\lambda_k^{p-1}}: C^{\infty}_c(c,+\infty)\oplus
C^{\infty}_c(c,+\infty)\longrightarrow L^2(c,+\infty)\oplus
L^2(c,+\infty)$$ given by (\ref{w1w2}). Moreover, for every
$i=1,2,3$, for every $k\in \mathbb N$ and for every $c>0$, we have
that $\sigma_{\rm ess}(D_{i\lambda_k})=\sigma_{\rm
ess}(D^F_{i\lambda_k})$.\par
\bigskip
Thus, much information about the essential spectrum of $\Delta_M$
can be recovered by the investigation of the essential spectra of
the selfadjoint operators $D^F_{1\lambda^p_k}$,
$D^F_{2\lambda^{p-1}_k}$ and $D^F_{3\lambda_k^{p-1}}$ for
arbitrarily large $c$. Since the Hodge $*$ operator isometrically
maps $p$-forms of type I onto $(N-p)$-forms of type II, it
suffices to consider the cases $i=1$ and $i=3$ . We remark that,
since the direct sums in (\ref{1delta}) and (\ref{3}) have an
infinite number of summands, for $i=1,3$ $$\sigma_{\rm
ess}(\Delta_{Mi})\supset \bigcup_k \sigma_{\rm
ess}(D^F_{i\lambda_k})$$ but we cannot argue that $$\sigma_{\rm
ess}(\Delta_{Mi})= \bigcup_k \sigma_{\rm ess}(D^F_{i\lambda_k}).$$
\section{The essential spectrum}
In the present section, we will compute the essential spectrum of
$\Delta_M$, under suitable assumptions on the asymptotic behaviour
of $f$ and $g$. Namely, if $$\tilde f(t):=f(t)-1, $$ $$\tilde
g(t):=g(t)- \sinh^2t;$$ we will suppose that for $t>>0$
\begin{equation}\label{asymestg}|\tilde g(t)|\leq \frac{C}{t},\quad
|\frac{\partial \tilde g}{\partial t}|\leq \frac{C}{t},\quad
|\frac{\partial^2\tilde g}{\partial t^2}|\leq
\frac{C}{t},\end{equation}
\begin{equation}\label{asymestf}|\tilde f(t)|\leq
\frac{C}{t}, \quad |\frac{\partial \tilde f}{\partial t}|\leq
\frac{C}{t},\quad|\frac{\partial^2\tilde f}{\partial t}|\leq
\frac{C}{t}.\end{equation}\par \bigskip For $i=1,2,3$ and for
every $k\in \mathbb N$, let $\Delta_{Mi}$ and $D_{i\lambda_k}$ be
defined as in Section 4.\par  First of all, we will determine the
essential spectrum of $\Delta_{M1}$. To this purpose, let us
recall some basic facts.
\begin{defn}(\cite{Reed-Simon}) Let $A$ be a selfadjoint operator on a
Hilbert space ${\mathcal H}$. An operator $C$ such that ${\mathcal
D}(A)\subset {\mathcal D}(C)$ is called relatively compact with
respect to $A$ if and only if $C(A+ i I)^{-1}$ is compact.
\end{defn} In terms of the Hilbert space ${\mathcal D}(A)$ endowed
with the norm $\|\phi\|_A$ given by $$\|\phi\|^2_A =
\|\phi\|_{\mathcal H}^2+ \|A\phi\|_{\mathcal H}^2, $$ $C$ is
relatively compact if and only if $C$ is compact from ${\mathcal
D}(A)$ with the norm $\|.\|_A$ to ${\mathcal H}$ with the norm
$\|.\|_{\mathcal H}$. Moreover, we recall the following Lemma (for
a proof see \cite{Reed-Simon}):
\begin{lem}\label{relcomp} Let $A$ be a selfadjoint operator on a Hilbert space
${\mathcal H}$, and let $C$ be a symmetric operator such that $C$ is a
relatively compact perturbation for $A^n$ for some positive integer $n$.
Suppose further that $B=A+C$ is selfadjoint on ${\mathcal D}(A)$, Then
$$ \sigma_{\rm ess}(A)= \sigma_{\rm ess}(B) .$$
 \end{lem}
Finally, we recall that, given a selfadjoint operator
$A$ on a Hilbert space ${\mathcal H}$, $\mu \in \sigma_{\rm
ess}(A)$ if and only if there exists a Weyl sequence
$\left\{w_n\right\}\subset {\mathcal D}(A)$ for $\mu$, that is, a
sequence $\left\{w_n\right\}\subset {\mathcal D}(A)$ with no
convergent subsequences in ${\mathcal H}$, bounded in ${\mathcal
H}$ and such that $$\lim_{n\rightarrow +\infty}(A-\mu) w_n =0
\quad\hbox{in ${\mathcal H}$}.$$
We are now in position to prove
our first result.
\begin{lem}\label{D1lambda} Let $M$ be endowed with a Riemannian metric
of type
(\ref{metric}), satisfying conditions (\ref{asymestg}),
(\ref{asymestf}); then for $0\leq p \leq N-1$,
 $$\sigma_{\rm
ess}(D_{1\lambda_k^p}^F)= \left[
\left(\frac{N-2p-1}{2}\right)^2,+\infty \right)$$
for every $k\in \mathbb N$.
\end{lem}
\begin{proof}
Let us consider the Friedrichs extension $D_{10}^F$ of the
operator  with constant coefficients  $$
D_{10}:C^{\infty}_c(c,+\infty)\longrightarrow L^2(c,+\infty)$$
\begin{equation}\label{nonpert1} D_{10}w= -
\frac{\partial^2}{\partial t^2}w+ \frac{(N-1-2p)^2}{4} w.
\end{equation} It is well-known that $\sigma_{\rm
ess}(D_{10}^F)=\left[ \left(\frac{N-2p-1}{2}\right)^2,+\infty
\right).$ We will show that $D^F_{1\lambda_k^p}-D^F_{10}$ is a
relatively compact perturbation of $(D^F_{10})^2$ for every $k \in
\mathbb N$. This, thanks to Lemma \ref{relcomp}, will give the
conclusion.
\par First of all, it is not difficult to see that for every $k\in
\mathbb N$, $${\mathcal D}((D^F_{10})^2)\subset {\mathcal
D}(D^F_{1\lambda_k^p}-D^F_{10});$$ indeed, comparing the domains
of $D_{1\lambda_k^p}^F$ and of $D_{10}^F$, we find that
$${\mathcal D}(D_{10}^F)={\mathcal D}(D^F_{1\lambda_k^p}).$$ We
still have to check that for every sequence $\left\{w_n
\right\}\subset {\mathcal D}((D^F_{10})^2)$ such that
\begin{equation}\label{bound1}\|w_n\|^2_{L^2}+
\|(D^F_{10})^2w_n\|^2_{L^2}\leq C \end{equation} there exists a
subsequence $\left\{ w_{n_l} \right\}$ such that
$\left\{(D^F_{1\lambda_k^p}- D^F_{10})w_{n_l}\right\}$ converges
in $L^2(c,+\infty)$.\par To this purpose, let us observe that
conditions (\ref{asymestg}) and (\ref{asymestf}) yield:
\begin{equation}\label{2inf11}\left( 1-\frac{1}{f} \right) \in
L^2(c,+\infty)\cap L^{\infty}(c, +\infty);
\end{equation}
\begin{equation}\label{2inf12} \frac{1}{f^2}\frac{\partial f}{\partial
t}\in
L^2(c,+\infty)\cap L^{\infty}(c, +\infty);
\end{equation}
\begin{equation}\label{2inf13}
W_1(t)\in L^2(c,+\infty)\cap L^{\infty}(c, +\infty),
\end{equation}
where
\begin{multline}W_1(t):=  \left\{
-\frac{7}{16} \frac{1}{f^3} \left(\frac{\partial f}{\partial
t}\right)^2+ \frac{1}{4} \frac{1}{f^2} \frac{\partial^2
f}{\partial t^2}+\right.\\ - \frac{1}{2}
\frac{1}{f^2}\frac{\partial f}{\partial t}\frac{(N-1-2p)}{4}
\frac{1}{g}\frac{\partial g}{\partial t} + \frac{1}{f}
\frac{(N-2p-1)}{4}
\frac{(N-2p-5)}{4}\frac{1}{g^2}\left(\frac{\partial g}{\partial
t}\right)^2\\ \left. +
\frac{1}{f}\frac{(N-2p-1)}{4}\frac{1}{g}\frac{\partial^2 g
}{\partial t^2} + \frac{\lambda}{g}\right\} -
\frac{(N-2p-1)^2}{4}.
\end{multline}
Moreover, by (\ref{bound1}), the sequence $\left\{w_n \right\}$ is
bounded in $W^{3,2}(c,+\infty)$, and the Sobolev embedding theorem
implies that $\left\{w_n\right\}$, $\left\{\frac{\partial
w_n}{\partial t}\right\}$, $\left\{\frac{\partial^2 w_n}{\partial
t^2}\right\}$ are bounded sequences in $L^{\infty}(c, +\infty)$.
\par  Now, for every $n,m\in \mathbb N$
\begin{multline}\|(D^F_{1\lambda_k^p}-D^F_{10})(w_n-w_m
)\|_{L^2(c,+\infty)}\leq \\
\|\left(1-\frac{1}{f}\right)\frac{\partial^2}{\partial
t^2}(w_n-w_m)\|_{L^2(c,+\infty)}+\|\frac{\partial
f}{\partial t}(w_n-w_m)\|_{L^2(c,+\infty)}\\ +
\|W_1(w_n-w_m)\|_{L^2(c,+\infty)}.
\end{multline}
Let us begin with the third summand. For any compact subset
$K\subset (c,+\infty)$ and for every $n,m\in \mathbb N$ $$ \|W_1
(w_n-w_m)\|^2_{L^2(c,+\infty)} \leq C\int_K(w_n-w_m)^2 \,ds +
C\int_{(c,+\infty)\setminus K}W_1^2(s)\,ds, $$ where $C$ is a
positive constant independent of $K$. Indeed,
$W_1^2\in L^{\infty}(c,+\infty)$ and $(w_n-w_m)^2$ is bounded in
$L^{\infty}(c,+\infty)$. \par Let us consider a sequence
$\left\{c_h\right\}\subset (c,+\infty)$ such that $c_h \rightarrow
+\infty$ as $h\rightarrow +\infty$ and for every $h\in \mathbb N$
$$C\int_{c_h}^{+\infty}W_1^2(s)\,ds < \frac{1}{h} .$$ For $h=1$,
thanks to the Rellich-Kondrachov theorem, there exists a
subsequence $\left\{w_{n(1)}\right\}$ such that
$\left\{(w_{n(1)})_{|(c,c_1)}\right\}$ converges in $L^2(c,c_1)$.
Hence, for every $\delta
>0$ there exists $\bar n(1)$ such that for every $n,m
>\bar n(1)$ $$\int_c^{+\infty}W_1^2(w_{n(1)}-w_{m(1)})^2\,ds<
\frac{\delta}{3} +1. $$ Analogously, for $h=2$ there exists a
subsequence $\left\{w_{n(2)}\right\}\subseteq \left\{
w_{n(1)}\right\}$ such that for every $\delta >0$ there exists
$\bar n(2)$ such that for every $n,m>\bar n(2)$
$$\int_c^{+\infty}W_1^2(w_{n(2)}-w_{m(2)})^2\,ds< \frac{\delta}{3}
+\frac{1}{2}. $$ Going on in this way, for every $h\in \mathbb N$
we can find a subsequence $\left\{w_{n(h)}\right\} \subseteq
\left\{ w_{n(h-1)}\right\}$ such that for every $\delta >0$ there
exists $\bar n(h)$ such that for every $n,m>\bar n(h)$
$$\int_c^{+\infty}W_1^2(w_{n(h)}-w_{m(h)})^2\,ds< \frac{\delta}{3}
+\frac{1}{h}. $$ Through a Cantor diagonal process, then, we can
find a subsequence $\left\{w_{n_l}\right\}\subseteq \left\{ w_n
\right\}$ such that $\left\{ W_1 w_{n_l}\right\}$ is a Cauchy
sequence in $L^2(c,+\infty)$.\par As for the estimates of the
other two summands, recalling (\ref{2inf11}) and (\ref{2inf12}),
since $\left\{\frac{\partial w_{n_l}}{\partial t}\right\}$ and
$\left\{\frac{\partial^2 w_{n_l}}{\partial t^2}\right\}$ are
bounded in $L^{\infty}(c,+\infty)$ and in $W^{1,2}(K)$ for any
compact set $K\subset (c,+\infty)$, we can apply  the same
procedure. As a consequence, we can extract a subsequence, again
denoted by $\left\{w_{n_l}\right\}$, such that
$\left\{(D^F_{1\lambda_k^p}-D^F_{10})w_{n_l}\right\}$ converges in
$L^2(c,+\infty)$. This yields the conclusion.
\end{proof}
As a consequence, $$\left[\left(\frac{N-2p-1}{2}\right)^2,
+\infty\right)\subset \sigma_{\rm ess}(\Delta_{M1}).$$ \par \bigskip On
the other hand, the following Lemma holds:
\begin{lem}\label{ess1} Let $M$ be endowed with a Riemannian
metric of type (\ref{metric}), such that $f(t)\rightarrow 1$ and
$g(t) \rightarrow \sinh^2 t$ as $t\rightarrow +\infty$; for $0\leq
p \leq N-1$, if
$\mu < (\frac{N-2p-1}{2})^2$, then $\mu \notin \sigma_{\rm ess}
(\Delta_{M1}).$
\end{lem}
\begin{proof}
First of all, for every $k\in \mathbb N$ let us write
$D'_{1\lambda_k^p}$ as $$D'_{1\lambda_k^p}w=
-\frac{\partial}{\partial t}\left(\frac{1}{f}\frac{\partial
w}{\partial t}\right)+\left(V_1(t)+\frac{\lambda_k^p}{g}\right)w,
$$ where
\begin{multline}V_1(t):=\left\{
-\frac{7}{16} \frac{1}{f^3} \left(\frac{\partial f}{\partial
t}\right)^2+ \frac{1}{4} \frac{1}{f^2} \frac{\partial^2
f}{\partial t^2}\right.\\ - \frac{1}{2}
\frac{1}{f^2}\frac{\partial f}{\partial t}\frac{(N-1-2p)}{4}
\frac{1}{g}\frac{\partial g}{\partial t} + \frac{1}{f}
\frac{(N-2p-1)}{4}
\frac{(N-2p-5)}{4}\frac{1}{g^2}\left(\frac{\partial g}{\partial
t}\right)^2\\ \left. +
\frac{1}{f}\frac{(N-2p-1)}{4}\frac{1}{g}\frac{\partial^2 g
}{\partial t^2} \right\}.
\end{multline}
Now, let $\mu < \left(\frac{N-1-2p}{2} \right)^2$. Since
 for every $k\in \mathbb N$ the essential spectrum
of $D^F_{1\lambda_k^p}$ does not depend on $c$ and since $V_1(t)$
converges to $\left(\frac{N-2p-1}{2}\right)^2>\mu$ as
$t\rightarrow +\infty$, we can choose $c>0$ such that for every
$t>c$ $$V_1(t)-\mu
>C$$ for some positive constant $C>0$. \par If $\mu \in \sigma_{\rm
ess}(\Delta_{M1})=\sigma_{\rm ess}(\Delta_{M1}^F)$, there exists a
Weyl sequence for $\mu$, that is a sequence $\left\{\omega_k
\right\}\subset {\mathcal D}(\Delta_{M1}^F)$ such that $$\langle
\omega_k,\omega_k\rangle_{L^2_p(M)}\leq C ,$$ $$\lim_{k
\rightarrow +\infty}(\Delta_{M1}^F\omega_k -\mu \omega_k)=0 ,$$
from which it is not possible to extract any subsequence
converging in $L^2_p(M)$. Moreover, we can suppose that
$$\omega_k= h_k(t) \tau_{1k},$$ where $\tau_{1k}$ is a coclosed
normalized $p$-eigenform of $\Delta_{\Sphere^{N-1}}$ corresponding
to $\lambda_k^p$ and $\lambda_k^p \rightarrow +\infty$ as
$k\rightarrow +\infty$. Hence, via unitary equivalence, there
exists a sequence $\left\{w_k\right\}\subset {\mathcal
D}(D^F_{1\lambda_k^p}) $ such that $$\|w_k\|_{L^2(c,+\infty)}\leq
C $$
\begin{equation}\label{Weyl1}\lim_{k\rightarrow
+\infty}\|D^F_{1\lambda_k^p}w_k -\mu w_k\|_{L^2(c,+\infty)}=0,
\end{equation}
from which we cannot extract any $L^2$-converging subsequence.
Then $$\langle D^F_{1\lambda_k^p}w_k -\mu w_k ,w_k
\rangle_{L^2(c,+\infty)}\longrightarrow 0 $$ as $k\rightarrow
+\infty$, and, since for every $k\in \mathbb N$ $${\mathcal
D}(D^F_{1\lambda_k^p})\subset W^{1,2}_0(c,+\infty),$$ we get
\begin{multline}\int_{c}^{+\infty}\frac{1}{f(s)}
\left(\frac{\partial w_k}{\partial s}\right)^2(s)\,ds +
\int_{c}^{+\infty} [V_1(s)-\mu]w_k^2(s)\,ds \\+
\int_{c}^{+\infty}\frac{\lambda_k^p}{g(s)}w_k^2(s)\,ds
\longrightarrow 0\end{multline} as $k\rightarrow +\infty$. Since
all the terms are positive, we have
$$\int_{c}^{+\infty}[V_1(s)-\mu]w_k^2(s)\,ds\longrightarrow 0  $$
as $k\rightarrow +\infty$, whence $$\int_{c}^{+\infty}w_k^2(s)\,ds
\longrightarrow 0$$ as $k \rightarrow +\infty$, because
$$\int_{c}^{+\infty}w_k^2(s)\,ds \leq
\frac{1}{C}\int_{c}^{+\infty}[V_1(s)-\mu]w_k^2(s)\,ds.$$ This
yields a contradiction. Hence, if $\mu \leq
\left(\frac{N-2p-1}{2}\right)^2$, $\mu\notin \sigma_{\rm
ess}(\Delta_{M1})$. \end{proof} As a consequence
\begin{prop}\label{spess1} Let $M$ be endowed with a Riemannian metric of
type
(\ref{metric}), satisfying conditions (\ref{asymestg}),
(\ref{asymestf}); then, for $0\leq p \leq N-1$, $$\sigma_{\rm
ess}(\Delta_{M1})=\left[\left(\frac{N-2p-1}{2}\right)^2,+\infty\right).$$
\end{prop}
By duality,
\begin{prop}\label{spess2} Let $M$ be endowed with a Riemannian metric of
type
(\ref{metric}), satisfying conditions (\ref{asymestg}),
(\ref{asymestf}); then, for $1\leq p \leq N$,  $$\sigma_{\rm
ess}(\Delta_{M2})=\left[\left(\frac{N-2p+1}{2}\right)^2,+\infty\right).$$
\end{prop}
We still have to determine the essential spectrum of
$\Delta_{M3}$ for $1\leq p \leq N-1$. First of all, we compute the
essential spectrum of
$D^F_{3\lambda_k^{p-1}}$ for every $k\in \mathbb N$:
\begin{lem}\label{D3lambda} Let $M$ be endowed with a Riemannian
metric of type (\ref{metric}), satisfying conditions
(\ref{asymestg}) and (\ref{asymestf}); then, for $1\leq p \leq N-1$,
$$\sigma_{\rm
ess}(D^F_{3\lambda_k^{p-1}})=
\left[\min\left\{\left(\frac{N-2p-1}{2}\right)^2, \left(
\frac{N-2p+1}{2} \right)^2 \right\},+\infty \right) $$ for every
$k\in \mathbb N$.
\end{lem}
\begin{proof} Let us consider the Friedrichs extension
$D_{30}^F$ of the  operator
$$D_{30}:C^{\infty}_c(c,+\infty)\oplus
C^{\infty}_c(c,+\infty)\longrightarrow L^2(c,+\infty)\oplus
L^2(c,+\infty)$$
\begin{multline}D_{30} (w_1 \oplus w_2):=
\left(-\frac{\partial^2 w_1}{\partial t^2}+
\left(\frac{N-2p-1}{2}\right)^2 w_1 \right)\\ \oplus \left(
-\frac{\partial^2 w_2}{\partial
t^2}+\left(\frac{N-2p+1}{2}\right)^2 w_2 \right) \end{multline} It
is not difficult to see that $$\sigma_{\rm ess}(D^F_{30})=
\left[\min\left\{\left(\frac{N-2p-1}{2}\right)^2, \left(
\frac{N-2p+1}{2} \right)^2 \right\},+\infty \right). $$ As in
Lemma \ref{D1lambda}, we will show that
$D^F_{3\lambda_k^{p-1}}-D^F_{30}$ is a relatively compact
perturbation of $(D^F_{30})^2$. First of all, ${\mathcal
D}((D^F_{30})^2)\subset {\mathcal D}(D_{30}^F)={\mathcal
D}(D^F_{3\lambda_k^{p-1}}-D^F_{30})$; indeed, an explicit
comparison of the domains shows that ${\mathcal
D}(D_{30}^F)={\mathcal D}(D^F_{3\lambda_k^{p-1}})$ for every $k\in
\mathbb N$. \par We still have to check that for every sequence
 $$\left\{w_{1n}\oplus w_{2n} \right\}\subset
{\mathcal D}((D^F_{30})^2)$$ such that
\begin{equation}\label{bounded3}\|w_{1n}\oplus w_{2n}\|_{L^2 \oplus L^2}+
\|(D^F_{30})^2(w_{1n}\oplus w_{2n})\|_{L^2\oplus L^2}\leq C
,\end{equation} there exists a subsequence $\left\{w_{1n_l}\oplus
w_{2n_l}\right\}$ such that
$$(D^F_{3\lambda_k^{p-1}}-D^F_{30})(w_{1n_l}\oplus w_{2n_l})$$
converges in $L^2(c,+\infty)\oplus L^2(c,+\infty)$.  \par Now,
(\ref{bounded3}) implies that $\left\{ w_{in}\right\}$
is bounded in $W^{3,2}(c,+\infty)$ for $i=1,2$; hence
$\left\{w_{in}\right\}$,
$\left\{\frac{\partial w_{in}}{\partial t}\right\}$ and
$\left\{\frac{\partial^2 w_{in}}{\partial t^2}\right\}$ are
bounded in $L^{\infty}(c,+\infty)$ and in $W^{1,2}(K)$ for every
compact subset $K\subset (c,+\infty)$. For every $n,m\in
\mathbb N$ $$\|(D^F_{3\lambda_k^{p-1}}-D^F_{30})((w_{1n}-w_{1m})
\oplus(w_{2n}-w_{2m}))\|_{L^2(c,+\infty)\oplus
L^2(c,+\infty)}\leq$$ $$\leq
\|(D^F_{1\lambda_k^{p-1}}-D^F_{10})(w_{1n}-w_{1m})\|_{L^2(c,+\infty)}
+$$ $$+
\|(D^F_{2\lambda_k^{p-1}}-D^F_{20})(w_{2n}-w_{2m})\|_{L^2(c,+\infty)}+$$
$$ + \|V_{3\lambda_k^{p-1}}(w_{1n}-w_{1m})\|_{L^2(c,+\infty)} +
\|V_{3\lambda_k^{p-1}}(w_{2n}-w_{2m})\|_{L^2(c,+\infty)},$$ where
$$ V_{3\lambda_k^{p-1}}(t):=
g(t)^{-\frac{3}{2}}f(t)^{-\frac{1}{2}}\frac{\partial g}{\partial
t}\sqrt{\lambda_k^{p-1}}.$$ The first two terms can be estimated
as in Lemma \ref{D1lambda}; as for the last two terms, since under
conditions (\ref{asymestg}) and (\ref{asymestf})
$$V_{3\lambda_k^{p-1}} \in L^2(c,+\infty)\cap L^\infty(c,+\infty)
,$$ following the argument of Lemma \ref{D1lambda} we get the
conclusion.
\end{proof}
We still have to check whether the essential spectrum of
$\Delta_{M3}$ can contain any other $\mu \in \Real$. The
techniques of Lemma \ref{ess1} can not be applied in this case,
because $D_{3\lambda_k^{p-1}}$ is a coupled system of differential
operators. Hence different techniques are needed. We have
\begin{lem}\label{ess3} Let $M$ be endowed with a Riemannian metric of type (\ref{metric}),
 satisfying conditions (\ref{asymestg}) and (\ref{asymestf}). For
$1\leq p \leq N-1$, if $0<\mu < \min \left\{
(\frac{N-2p-1}{2})^2,(\frac{N-2p+1}{2})^2\right\}$, then $\mu \notin
\sigma_{\rm ess} (\Delta_{M3}).$
\end{lem}
\begin{proof} We already know from Lemma \ref{D3lambda} that for every
$k \in \mathbb N$, $$\sigma_{\rm ess}(D^F_{3\lambda_k^{p-1}})=
\left[\min \left\{
\left(\frac{N-2p-1}{2}\right)^2,\left(\frac{N-2p+1}{2}\right)^2
\right\},+\infty\right).$$ As a consequence, given a positive $\mu
< \min \left\{ (\frac{N-2p-1}{2})^2,(\frac{N-2p+1}{2})^2\right\}$,
$\mu$ belongs to the essential spectrum of $\Delta_{M3}$ if and
only if there exist a sequence
 $\left\{\mu_k\right\}$ of eigenvalues of $\Delta_M$ and a corresponding
sequence $\left\{\Phi_k \right\}$ of $p$-forms of type III such
that $$ \mu_k \longrightarrow \mu \quad \hbox{as $k\longrightarrow
+\infty$}$$ and for every $k \in \mathbb N$ $$\Delta_M \Phi_k -\mu_k
\Phi_k=0.$$
Since $\mu >0$, we can suppose, up to the choice of a subsequence,
that either for every $d_M \Phi_k \not=0$ for every
$k\in \mathbb N$ or $\delta_M \Phi_k \not=0$ for every $k\in \mathbb N$.
Let us suppose to be in the first case. In view of (\ref{dM1}), $d_M
\Phi_k$ is a $(p+1)$-form of type II; moreover,
$$ \|d_M \Phi_k \|_{L^2(M)}<C \quad \hbox{for every $k\in \mathbb
N$}$$
$$\Delta_M d_M \Phi_k -\mu_k d_M\Phi_k =0  \quad \hbox{for every $k\in
\mathbb
N$}, $$
and $$ \mu_k \longrightarrow \mu \quad \hbox{as $k
\rightarrow +\infty$ }.$$ Hence, $\mu \in \sigma_{\rm
ess}(\Delta_{M2})$, and, thanks to Proposition \ref{spess2}, $$\mu >
\left(
\frac{N-2(p+1)+1}{2}\right)^2=\left(\frac{N-2p-1}{2}\right)^2,$$
in contradiction with our hypothesis.\par If on the contrary we
are in the second case, in view of (\ref{deltaM1}), $\delta_M
\Phi_k$ is a $(p-1)$-form of type I; moreover,
$$\|\delta_M \Phi_k \|_{L^2(M)}<C \quad \hbox{for every $k\in
\mathbb N$},$$
$$\Delta_M \delta_M
\Phi_k -\mu_k \delta_M \Phi_k =0 \quad \hbox{for every $k\in \mathbb
N$}, $$
and $$ \mu_k \longrightarrow
\mu \quad \hbox{as $k \rightarrow +\infty$ }.$$ Hence, $\mu \in
\sigma_{\rm ess}(\Delta_{M1})$, and by Proposition \ref{spess1} $$\mu
> \left(
\frac{N-2(p-1)-1}{2}\right)^2=\left(\frac{N-2p+1}{2}\right)^2,$$
in contradiction with our hypothesis.\par This yields the
conclusion.
\end{proof}
Hence,
\begin{prop}\label{spess3} Let $M$ be endowed with a Riemannian metric of
type (\ref{metric}), satisfying conditions (\ref{asymestg}),
(\ref{asymestf}); then, for $1\leq p \leq N-1$, $$\sigma_{\rm
ess}(\Delta_{M3})\setminus
\left\{0\right\}=\left[\min\left\{\left(\frac{N-2p-1}{2}\right)^2
,\left(\frac{N-2p+1}{2}\right)^2 \right\}, +\infty \right). $$
\end{prop}
\begin{rem}\label{zeroacc} By the same argument of Lemma \ref{ess3} it is
possible to show that if there exist a sequence $\left\{\mu_k\right\}$ of
positive eigenvalues of $\Delta_M$ and a corresponding sequence
$\left\{\Phi_k\right\}$ of $p$-forms of type III such that
$$ \mu_k \longrightarrow 0 \quad \hbox{as $k\longrightarrow +\infty$}$$
and for every $k\in \mathbb N$
$$ \Delta_M \Phi_k -\mu_k \Phi_k =0,$$
then $0\notin \sigma_{\rm ess}(\Delta_{M3})$.
\end{rem}
Recalling the results of section 3, finally we can completely
determine the essential spectrum of $\Delta_M$:
\begin{thm}\label{main} Let $M$ be endowed with a Riemannian metric of
type
(\ref{metric}) satisfying condition (\ref{condzero}) and conditions
(\ref{asymestg}), (\ref{asymestf}). Then, if $p\not= \frac{N}{2}$,
$$\sigma_{\rm ess}(\Delta_M)= \left[\min \left\{
\left(\frac{N-2p-1}{2}\right)^2,\left(\frac{N-2p+1}{2}\right)^2\right\},
+\infty\right) $$
whilst if $p=\frac{N}{2}$
$$\sigma_{\rm ess}(\Delta_M)= \left\{0\right\}\cup \left[\frac{1}{4},
+\infty\right)
.$$
\end{thm}
\begin{proof} Thanks to Propositions \ref{spess1}, \ref{spess2},
\ref{spess3} we have that
$$ \sigma_{\rm ess}(\Delta_M) \setminus \left\{0\right\}=
\left[\min\left\{\left(\frac{N-2p-1}{2}\right)^2,\left(
\frac{N-2p+1}{2}\right)^2\right\},+\infty\right).$$
Moreover, in view of Remark \ref{zeroacc}, $0$ can belong to the essential
spectrum of $\Delta_M$ if and only if it is an eigenvalue of $\Delta_M$ of
infinite multiplicity. In view of Theorem \ref{zeroess}, this happens only
if $p=\frac{N}{2}$. The conclusion follows.
\end{proof}


\end{document}